\documentclass[11pt]{article}
\usepackage[dvips]{graphics,epsfig}
\usepackage[latin1]{inputenc}
\usepackage{latexsym,enumerate}
\usepackage{amsmath,amssymb}
\usepackage{pxfonts}

\def\cqfd{$\sqcup\!\!\!\!\sqcap$}

\newtheorem{theorem}{Theorem}

\newtheorem{lemma}{Lemma}
\newtheorem{corollary}{Corollary}

\begin{document}

\centerline{\bf BIAS-REDUCED ESTIMATORS OF THE WEIBULL TAIL-COEFFICIENT}
\vskip4ex
\large  Jean Diebolt$^{(1)}$, Laurent Gardes$^{(2)}$, St\'ephane Girard$^{(3)}$ and Armelle Guillou$^{(4)}$

\vskip4ex
\noindent
{\it $^{(1)}$ CNRS, Universit\'e de Marne-la-Vall\'ee\\
\'Equipe d'Analyse et de Math\'ematiques Appliqu\'ees\\
5, boulevard Descartes, Batiment Copernic\\
Champs-sur-Marne\\
77454 Marne-la-Vall\'ee Cedex 2}

\vskip2ex
\noindent
{\it $^{(2)}$
Universit\'e Grenoble 2,\\
LabSAD, 1251 Avenue centrale \\
B.P. 47, 38040 Grenoble Cedex 9}

\vskip2ex
\noindent
{\it $^{(3)}$
Universit\'e Grenoble 1,\\
LMC-IMAG, 51 rue des Math\'ematiques \\
B.P. 53, 38041 Grenoble Cedex 9}

\vskip2ex
\noindent
{\it $^{(4)}$ Universit\'e Paris VI\\
Laboratoire de Statistique Th\'eorique et Appliqu\'ee \\
Bo\^ \i te 158\\
175 rue du Chevaleret \\
75013 Paris}

\vskip4ex
\noindent
{\bf Abstract.} {\it In this paper, we consider the problem of the estimation of a Weibull tail-coefficient~$\theta$. In particular,
we propose a regression model, from which we derive a bias-reduced estimator of~$\theta$. This estimator is based on a least-squares approach. The asymptotic normality of this estimator is also established. A small simulation study is provided in order to prove its efficiency.}


\vskip8ex

\section{Introduction}\label{sec.un}

Let $X_1, ..., X_n$ be a sequence of independent and identically distributed random variables with distribution function $F$, and let $X_{1,n} \leq ... \leq X_{n,n}$ denote the order statistics associated to this sample.

In the present paper, we address the problem of estimating the Weibull tail-coefficient~$\theta>0$ defined as
\begin{eqnarray}
1-F(x) = \exp(-H(x)) \mbox{ with } H^{-1}(x) := \inf\{t: H(t) \geq x\}= x^{\theta} \ell(x),
\label{model}
\end{eqnarray}
where $\ell$ is a slowly varying function at infinity satisfying
\begin{equation}
\label{convl}
\frac{\ell(\lambda x) }{ \ell(x)} \longrightarrow 1, \, \mbox{as}\, x \to \infty, \, \mbox{for all}\, \lambda > 0.
\end{equation}
\cite{bookbeir} investigated this estimation problem and proposed the following estimator of~$\theta$:
\begin{equation}
\label{thetag}
\widetilde \theta_n(k_n) = {\displaystyle\frac{\sum_{i=1}^{k_n} \Bigr(\log X_{n-i+1,n}-\log X_{n-k_n+1,n}\Bigr)}  
{{\sum_{i=1}^{k_n} \Bigr(\log\log (n /i)-\log \log (n / k_n)\Bigr)}}},
\end{equation}
where $k_n$ is an intermediate sequence, i.e. a sequence such that $k_n \to \infty$ and ${k_n / n}\to 0$ as $n\to \infty$.

We refer to \cite{theta,Gar2,BeirBro} and \cite{BRO} for other 
propositions and to \cite{BeirBouq} for Local Asymptotic 
Normality (LAN) results.  
Estimator~(\ref{thetag}) is closed in spirit to the Hill estimator \cite{Hill} in the case of Pareto-type distributions. In~\cite{theta} , the asymptotic normality of $\widetilde \theta_n(k_n)$ is established under suitable assumptions. To prove such a result, a second-order condition is required in order to specify the bias-term. This assumption can be expressed in terms of the slowly varying function $\ell$ as follows:

{\bf Assumption $(R_{\ell}(b, \rho))$} {\it There exists a constant $\rho < 0$ and a rate function $b$ satisfying $b(x)\to 0$ as $x \to \infty$, such that for all $\varepsilon > 0$ and $1<A<\infty$, we have
$$\sup_{\lambda \in [1, A]} \Biggr\vert \frac{\log (\ell(\lambda x)/ \ell(x))}{ b(x) K_{\rho}(\lambda)}-1\Biggr\vert \leq \varepsilon,\quad \mbox{for $x$ sufficiently large,}$$
with $\displaystyle{K_{\rho}(\lambda) = \int_1^{\lambda} t^{\rho-1}dt.}$}

It can be shown that necessarily $\vert b\vert$ is regularly varying with 
index $\rho$ \cite{GEL}.
Moreover, we focus on the case where the convergence~(\ref{convl}) is slow, and thus when the
bias term in $\widetilde \theta_n(k_n)$ is large. This situation is described by the following
assumption:
\begin{equation}
\label{hypob}
x|b(x)|\to\infty \mbox{ as } x\to\infty.
\end{equation}
Let us note that this condition implies $\rho\geq -1$. Gamma and Gaussian distributions fulfill~(\ref{hypob}),
whereas Weibull distributions do not (see Table~\ref{tabex}) since, in this case, the bias term vanishes.

Using this framework, we will establish rigorously in Section~\ref{sec.deux} the following approximation for the log-spacings of upper order statistics:
\begin{eqnarray}
Z_j &:= &j\, \log(n/j) \Bigr(\log X_{n-j+1,n}-\log X_{n-j,n}\Bigr) \nonumber\\
&\approx& \Biggr(\theta + b\Bigr(\log (n / k_n)\Bigr)\Biggr(\frac{\log (n / j)}{\log(n / k_n)}\Biggr)^{\rho}\Biggr)f_j, 
\label{modele}
\end{eqnarray}
for $1\leq j \leq k_n$,
where $(f_1, ..., f_{k_n})$ is a vector of independent and standard exponentially distributed random variables.
This exponential regression model is similar to the ones proposed by 
\cite{BEIR,BDGS} and \cite{Feu} in the case of 
Pareto-type distributions. 
Ignoring $b(\log (n / k_n))\Bigr(\frac{\log (n / j)}{ \log(n / k_n)}\Bigr)^{\rho}$ in (\ref{modele}) leads to the maximum likelihood estimator
$$\check{\theta}_n(k_n) = \frac{1}{ k_n} \sum_{j=1}^{k_n} Z_j,$$ 
which turns out to be an alternative estimator of  $\widetilde \theta_n(k_n)$.
The full model (\ref{modele}) allows us to generate bias-corrected estimates 
$\widehat \theta_n(k_n)$ for $\theta$ through maximum likelihood estimation 
of~$\theta$, $b(\log {n / k_n})$ and $\rho$ for each $1\leq k_n \leq n-1.$ An 
alternative to this approach consists in using a canonical choice for $\rho$ 
and to estimate the two other parameters by a least-squares method (LS). For 
the canonical choice of $\rho$, we can use for instance the value -1, which is 
the same as the one proposed by \cite{Feu} for the regression 
model in the case of Pareto-type distributions. The asymptotic normality of 
the resulting LS-estimator is established in Section~\ref{sec.trois}.
An adaptive
selection method for $k_n$ in $\check\theta_n(k_n)$ is also derived.
In order to illustrate the usefulness of these results, we provide a 
simulation study in Section~\ref{sec.quatre} as well as an application to a real data set
in Section~\ref{sec.cinq}.
The proofs of our results are postponed to Section~\ref{sec.sept}.


\section{Exponential regression model} \label{sec.deux}

In this section, we formalize (\ref{modele}). First, remark that
$$F^{-1}(x) =[-\log(1-x)]^{\theta} \ell(-\log(1-x)).$$
Since $X_{n-j+1,n}  \stackrel{d}{=} F^{-1}(U_{n-j+1,n}), 1\leq j \leq k_n,$ where $U_{j,n}$ denotes the $j$-th order statistic of a uniform sample of size $n$, we have
$$X_{n-j+1,n} \stackrel{d}{=} \left[-\log(1-U_{n-j+1,n})\right]^{\theta}
\ell\left(-\log(1-U_{n-j+1,n})\right)$$
which implies that
$$\log X_{n-j+1,n} \stackrel{d}{=} \theta \log\left[-\log(1-U_{n-j+1,n})\right]
+ \log\left[\ell\left(-\log(1-U_{n-j+1,n})\right)\right].$$
Moreover, considering the order statistics from an independent standard exponential sample, $E_{n-j+1,n} \stackrel{d}{=} - \log(1-U_{n-j+1,n})$. Therefore
\begin{eqnarray*}
\log X_{n-j+1,n} &\stackrel{d}{=}& \theta \log(E_{n-j+1,n}) + \log\left[\ell(E_{n-j+1,n})\right]\\ &=:& A_{n}(j)+B_{n}(j).
\end{eqnarray*}
Recall that $Z_j=j \log(n/ j) \left(\log X_{n-j+1,n}-\log X_{n-j,n}\right)$,
$1\leq j \leq k_n$.  
Then, our basic result now reads as follows.
\begin{theorem}
\label{thun}
Suppose (\ref{model}) holds together with $(R_{\ell}(b, \rho))$ and (\ref{hypob}). Then, if $k_n \to \infty$ and ${\log k_n / \log n} \to 0$, we have
\begin{eqnarray}
\label{regression}
\sup_{1\leq j \leq k_n} \left\vert Z_j
- \left(\theta + b(\log (n / k_n))\left(\frac{\log (n / j)}{ \log(n / k_n)}\right)^{\rho}\right)f_j  \right\vert 
 = o_{\mathbb P}\Bigr(b(\log(n / k_n))\Bigr),
\end{eqnarray}
where $(f_1, ..., f_{k_n})$ is a vector of independent and standard exponentially distributed random variables.
\end{theorem}
\noindent The proof of this theorem is based on the following two lemmas:
\begin{lemma}
\label{lemAn}
Suppose (\ref{model}) holds together with $(R_{\ell}(b, \rho))$ and (\ref{hypob}). Then, if $k_n \to \infty$ and
${k_n/n} \to 0$, we have
$$
\sup_{1\leq j \leq k_n} \left\vert j \log(n/ j) \left[A_{n}(j)-A_{n}(j+1)\right]- \theta \, f_j  \right\vert = o_{\mathbb P}\left(b(\log(n / k_n))\right),
$$
\end{lemma}
\noindent and
\begin{lemma}
\label{lemBn}
Suppose (\ref{model}) holds together with $(R_{\ell}(b, \rho))$. 
Then, if $k_n \to \infty$ and ${\log k_n / \log n} \to 0$, we have
\begin{eqnarray*}
&&\sup_{1\leq j \leq k_n} \left\vert j \log(n/ j) \left[B_{n}(j)-B_{n}(j+1)\right]- b(\log (n / k_n))\left(\frac{\log (n / j)}{ \log(n / k_n)}\right)^{\rho}f_j  \right\vert \nonumber\\
&&= o_{\mathbb P}\left(b(\log(n / k_n))\right).
\end{eqnarray*}
\end{lemma}
\noindent The proof of these lemmas is postponed to Section~\ref{sec.sept}. 
\begin{corollary}
\label{Coro}
Under the assumptions of Theorem~\ref{thun}, we also have
\begin{eqnarray*}
\label{regression2}
\sup_{1\leq j \leq k_n} \left\vert Z_j- \left(\theta + b(\log (n / k_n))\left(\frac{\log (n /j)}{ \log(n / k_n)}\right)^{-1}\right)f_j  \right\vert 
 = o_{\mathbb P}\left(b(\log(n / k_n))\right),
\end{eqnarray*}
where $(f_1, ..., f_{k_n})$ is a vector of independent and standard exponentially distributed random variables.
\end{corollary}
\noindent This implies that one can plug the canonical choice $\rho=-1$
in the regression model~(\ref{regression}) without perturbing
the approximation. 
From model~(\ref{regression}) we can easily deduce the asymptotic normality of the estimator $\check{\theta}_n(k_n)$, given in the next theorem:
\begin{theorem}
\label{thdeux}
Suppose (\ref{model}) holds together with $(R_{\ell}(b, \rho))$ and (\ref{hypob}). Then, if $k_n \to \infty$, $\sqrt{k_n} b(\log ({n / k_n})) \to \lambda \in \mathbb R$ and, if $\lambda=0$, ${\log k_n / \log n} \to 0$, we have
$$\sqrt{k_n}\, \Biggr(\check{\theta}_n(k_n) - \theta-b(\log (n / k_n)) \frac{1}{ k_n} \sum_{j=1}^{k_n} \Biggr(\frac{\log (n / j) }{ \log (n / k_n)}\Biggr)^{\rho}\Biggr) \stackrel{d}{\longrightarrow} {\cal N}(0, \theta^2).$$
\end{theorem}
\noindent The Asymptotic Mean Squared Error (AMSE) associated to $\check{\theta}_n(k_n)$
is thus given by:
\begin{equation}
\label{AMSE}
AMSE(\check{\theta}_n(k_n))= \frac{\theta^2}{k_n}+\left( b(\log(n/k_n)) \frac{1}{k_n}
\sum_{j=1}^{k_n} \left( \frac{\log(n/j)}{\log(n/k_n)} \right)^{\rho}\right)^2. 
\end{equation}
This model (\ref{regression}) now plays the central role in the remainder of 
this paper. First, it allows us to generate bias-corrected estimates of 
$\theta$.  Second, it leads to the number of upper
order statistics $k_n$ to be used in $\check\theta_n(k_n)$ by minimizing 
the AMSE given by~(\ref{AMSE}) after replacing $\theta$, $b$ and $\rho$
by estimators.
These two points are described in the next section.

\section{Bias-reduced estimates of $\theta$ and adaptive selection of $k_n$}
\label{sec.trois}

In order to reduce the bias of the estimator $\check \theta_n(k_n)$, we can either 
estimate simultaneously $\theta, b(\log{n/ k_n})$ and $\rho$ by a maximum 
likelihood method or estimate $\theta$ and $b$ by a least-squares approach 
after substituting a canonical choice for $\rho$. In fact, this second-order 
parameter is difficult to estimate in practice and we can easily check by 
simulations that fixing its value does not much influence the result. This 
problem has already been discussed in \cite{BEIR,BDGS} and \cite{Feu}
where similar observations have been made in the 
case of Pareto-type distributions. The canonical choice $\rho=-1$ is often 
used although other choices could be motivated performing a model selection.

In all the sequel, we will estimate $\theta$ and $b(\log(n / k_n))$ by a LS-method after substituting $\rho$ with the value $-1$. In that case, we find the following LS-estimators:
\begin{eqnarray*}
\begin{cases} \widehat \theta_n(k_n) = \overline{Z}_{k_n} - \widehat b(\log(n / k_n)) \overline x_{k_n}  & \cr
&\cr
\displaystyle{\widehat b(\log(n / k_n)) = \frac{\sum_{j=1}^{k_n} (x_j-\overline x_{k_n})Z_j }{ \sum_{j=1}^{k_n} (x_j-\overline x_{k_n})^2} }& \cr
\end{cases}
\end{eqnarray*}
\noindent
where $x_j = \Bigr(\frac{\log (n/ j) }{ \log (n/ k_n)}\Bigr)^{-1}$, $\overline x_{k_n} = \frac{1 }{ k_n} \sum_{j=1}^{k_n} x_j$ and $\overline Z_{k_n} = 
\frac{1}{ k_n} \sum_{j=1}^{k_n} Z_j$.  
\noindent
Our next goal is to establish, under suitable assumptions, the asymptotic normality of $\widehat \theta_n(k_n)$. This is done in the following theorem.
\begin{theorem}
\label{thtrois}
Suppose (\ref{model}) holds together with $(R_{\ell}(b, \rho))$ and (\ref{hypob}). Then, if $k_n\to\infty$  such that 
\begin{eqnarray}
\label{cond2}
&&\frac{\sqrt {k_n}}{ \log(n / k_n)}b(\log(n / k_n)) \to \Lambda \in {\mathbb R} \mbox{ and, if } \Lambda = 0,\\
\label{cond2bis}
&&\frac{\log^2 k_n }{ \log(n / k_n)}\to 0   \mbox{ and } \frac{\sqrt{k_n} }{ \log(n / k_n)} \to \infty,
\end{eqnarray}
\noindent
we have
$$\frac{\sqrt {k_n} }{ \log(n / k_n)} \Bigr(\widehat \theta_n(k_n) - \theta\Bigr) \stackrel{d}{\longrightarrow} {\cal N}(0, \theta^2).$$
\end{theorem}
\noindent Remark that the rate of convergence of $\check{\theta}_n(k_n)$ is the same as the one of $\widehat \theta_n(k_n)$ in the cases where both $\lambda$ and $\Lambda$ are not equal to 0. 
The proof of this theorem is postponed to Section~\ref{sec.sept}.

We can also take benefit of the estimation of $b(\log n/k_n)$ by
estimating the AMSE given in~(\ref{AMSE}) by:
$$
\widehat{AMSE}(\check{\theta}_n(k_n))= \frac{(\widehat \theta_n(k_n))^2}{k_n}+\left( \widehat b(\log(n/k_n)) \frac{1}{k_n}
\sum_{j=1}^{k_n} \left( \frac{\log(n/j)}{\log(n/k_n)} \right)^{-1}\right)^2. 
$$
Then, the intermediate sequence $k_n$ can be selected by minimizing
the previous quantity: 
$$
\hat k_n=\arg\min_{k_n} \widehat{AMSE}(\check{\theta}_n(k_n)).
$$
This adaptive procedure for selecting the number of upper order
statistics is in the same spirit as the one used by \cite{MB}
in the context of the extreme value index estimation.

In order to illustrate the usefulness of the bias reduction and
of the selection procedure, we provide a simulation study in the next section.


\section{A simulation study}
\label{sec.quatre}

First, the finite sample performances of the estimators ${\widehat{\theta}}_n(k_n)$, ${\widetilde{\theta}}_n(k_n)$
and ${\check{\theta}}_n(k_n)$ are investigated on~6 different distributions: 
$\Gamma(0.25,1)$, $\Gamma(4,1)$, $|{\mathcal N}|(0,1)$, 
${\mathcal W}(0.25,0.25)$, ${\mathcal W}(4,4)$ and 
${\mathcal D}(1,0.5)$, see the appendix for the definition
of the latter distribution.  
We limit ourselves to these three estimators, since it is shown in \cite{theta} that
${\widetilde{\theta}}_n(k_n)$ gives better results than the other approaches
\cite{BeirBro,BRO}.
In each case, $N=100$ samples $({\mathcal X}_{n,i})_{i=1,\dots,N}$
of size $n=500$ were simulated.
On each sample $({\mathcal X}_{n,i})$, the estimates $\widehat{\theta}_{n,i}(k_n)$,
$\widetilde{\theta}_{n,i}(k_n)$
and $\check{\theta}_{n,i}(k_n)$
were computed for $k_n=2,\dots,360$.
Finally, the Hill-type plots were built by drawing the points
$$
\left( k_n, \frac{1}{N}\sum_{i=1}^N \widehat{\theta}_{n,i}(k_n)\right), \ 
\left( k_n, \frac{1}{N}\sum_{i=1}^N \widetilde{\theta}_{n,i}(k_n)\right) \ 
{\rm{and}} \ \left( k_n, \frac{1}{N}\sum_{i=1}^N \check{\theta}_{n,i}(k_n)\right).
$$
We also present the associated MSE plots
obtained by plotting the points
\begin{eqnarray*}
&&\left( k_n, \frac{1}{N}\sum_{i=1}^N \left(\widehat{\theta}_{n,i}(k_n)-\theta\right)^2\right), \ 
\left( k_n, \frac{1}{N}\sum_{i=1}^N 
\left(\widetilde{\theta}_{n,i}(k_n)-\theta\right)^2\right) {\rm{ and }} \\
&&\left( k_n, \frac{1}{N}\sum_{i=1}^N 
\left(\check{\theta}_{n,i}(k_n)-\theta\right)^2\right).
\end{eqnarray*}

\begin{figure}
\begin{center}
\begin{minipage}{0.85\textwidth}
\includegraphics*[scale=0.5]{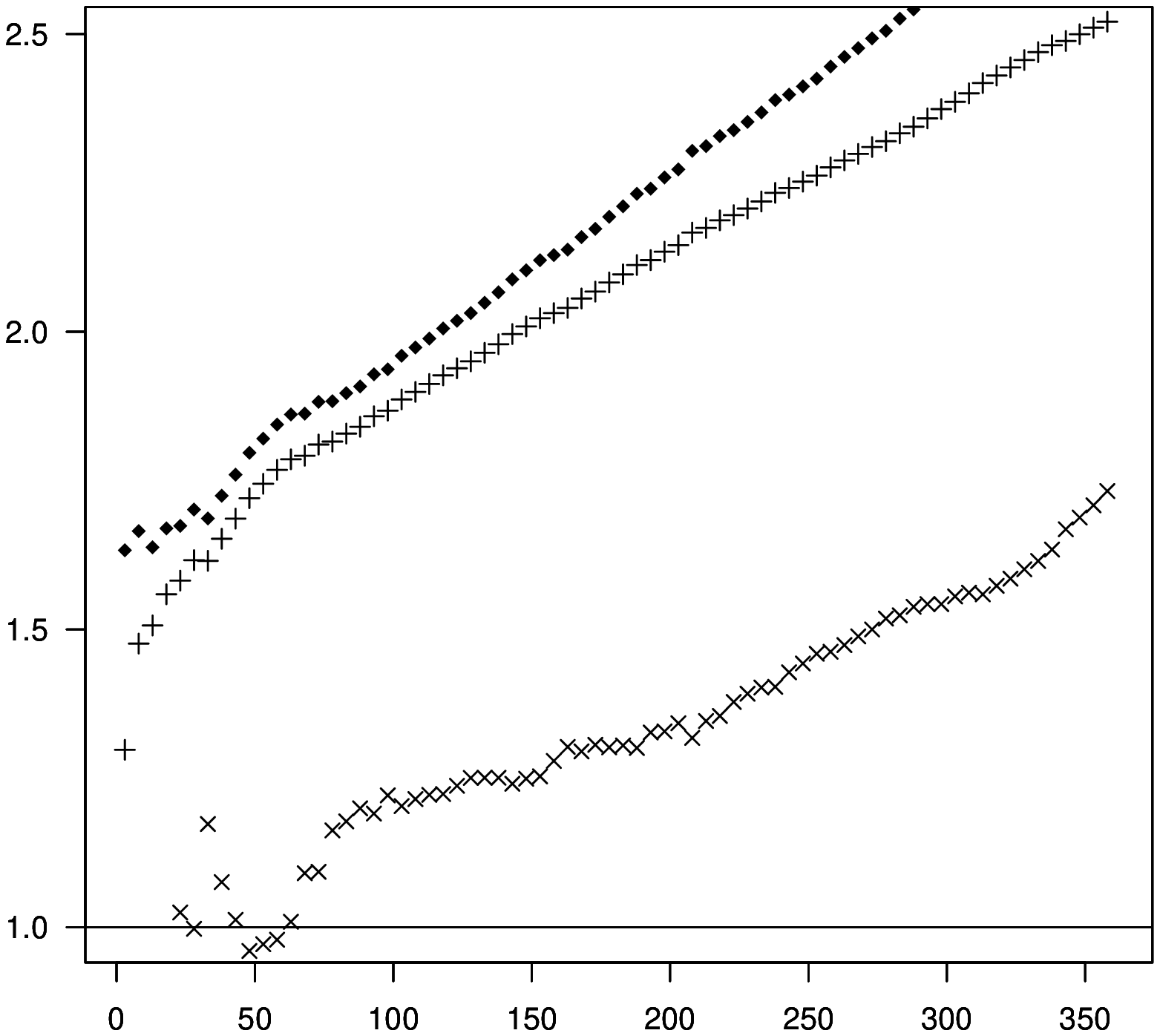}
\centerline{(a) Mean as a function of $k_n$}
\end{minipage}
\begin{minipage}{0.85\textwidth}
\includegraphics*[scale=0.5]{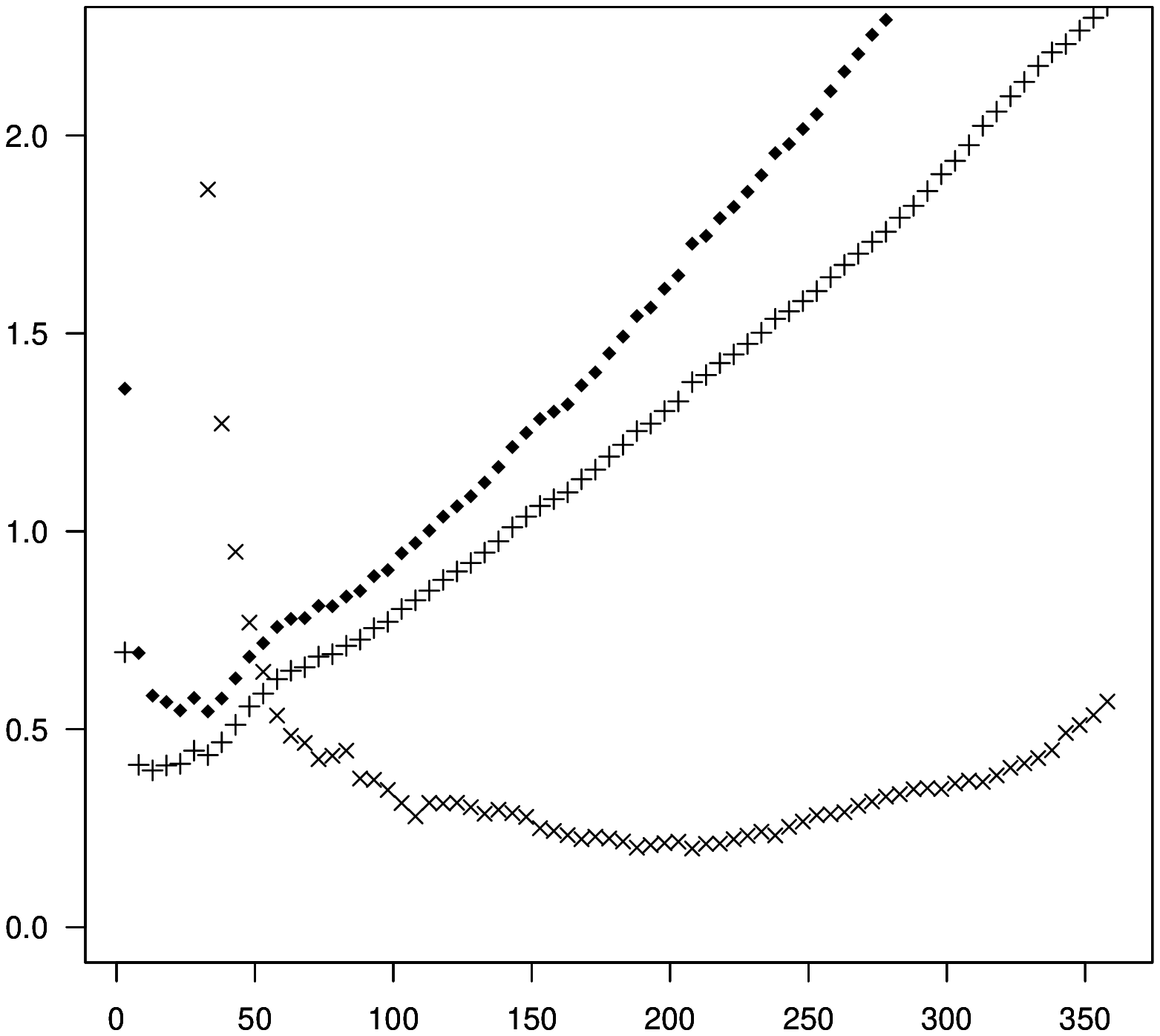}
\centerline{(b) Mean squared error as a function of $k_n$}
\end{minipage}
\end{center}
\caption{Comparison of estimates $\widehat{\theta}_n$ ($\times\times\times$), 
$\widetilde{\theta}_n$ ($\Diamondblack\Diamondblack\Diamondblack$) and $\check{\theta}_n$ ($+++$) for 
the $\Gamma(0.25,1)$ distribution. In (a), the straight line is the true value 
of $\theta$.}
\label{figgam025}
\end{figure}

\begin{figure}
\begin{center}
\begin{minipage}{0.85\textwidth}
\includegraphics*[scale=0.5]{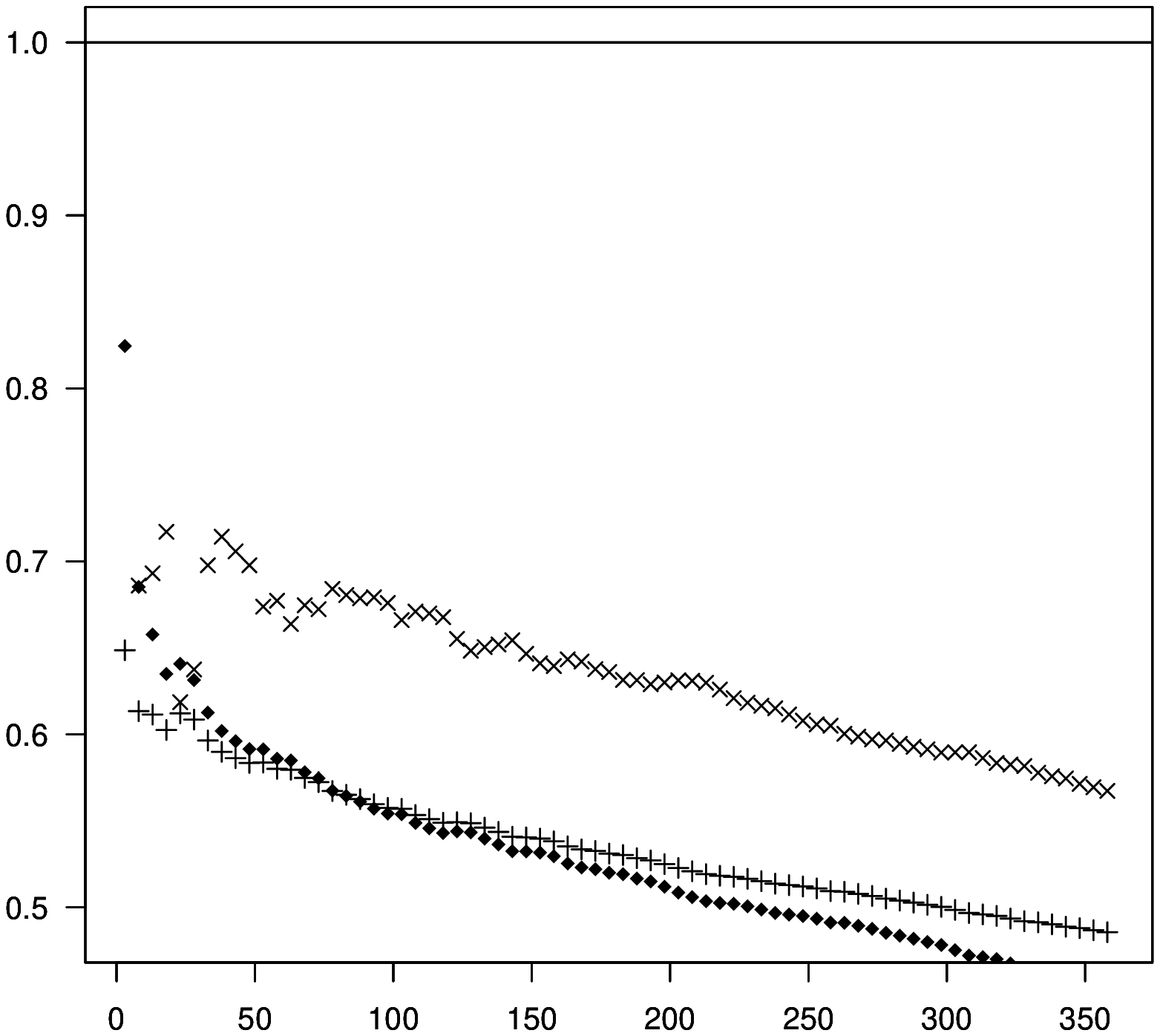}
\centerline{(a) Mean as a function of $k_n$}
\end{minipage}
\begin{minipage}{0.85\textwidth}
\includegraphics*[scale=0.5]{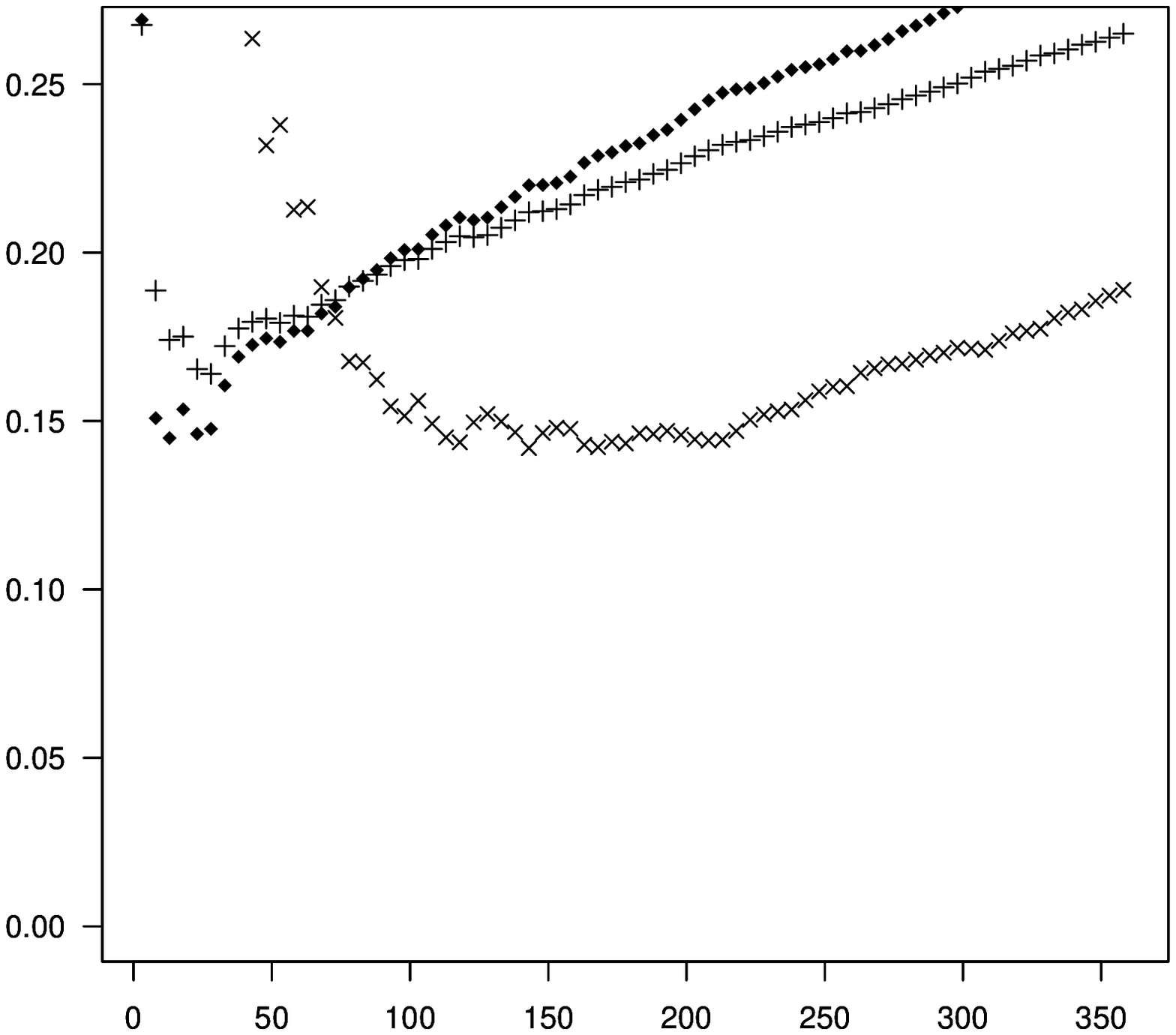}
\centerline{(b) Mean squared error as a function of $k_n$}
\end{minipage}
\end{center}
\caption{Comparison of estimates $\widehat{\theta}_n$ ($\times\times\times$), 
$\widetilde{\theta}_n$ ($\Diamondblack\Diamondblack\Diamondblack$) and $\check{\theta}_n$ ($+++$) for the $\Gamma(4,1)$ distribution. In (a), the straight line is the true value of $\theta$.}
\label{figgam4}
\end{figure}

\begin{figure}
\begin{center}
\begin{minipage}{0.85\textwidth}
\includegraphics*[scale=0.5]{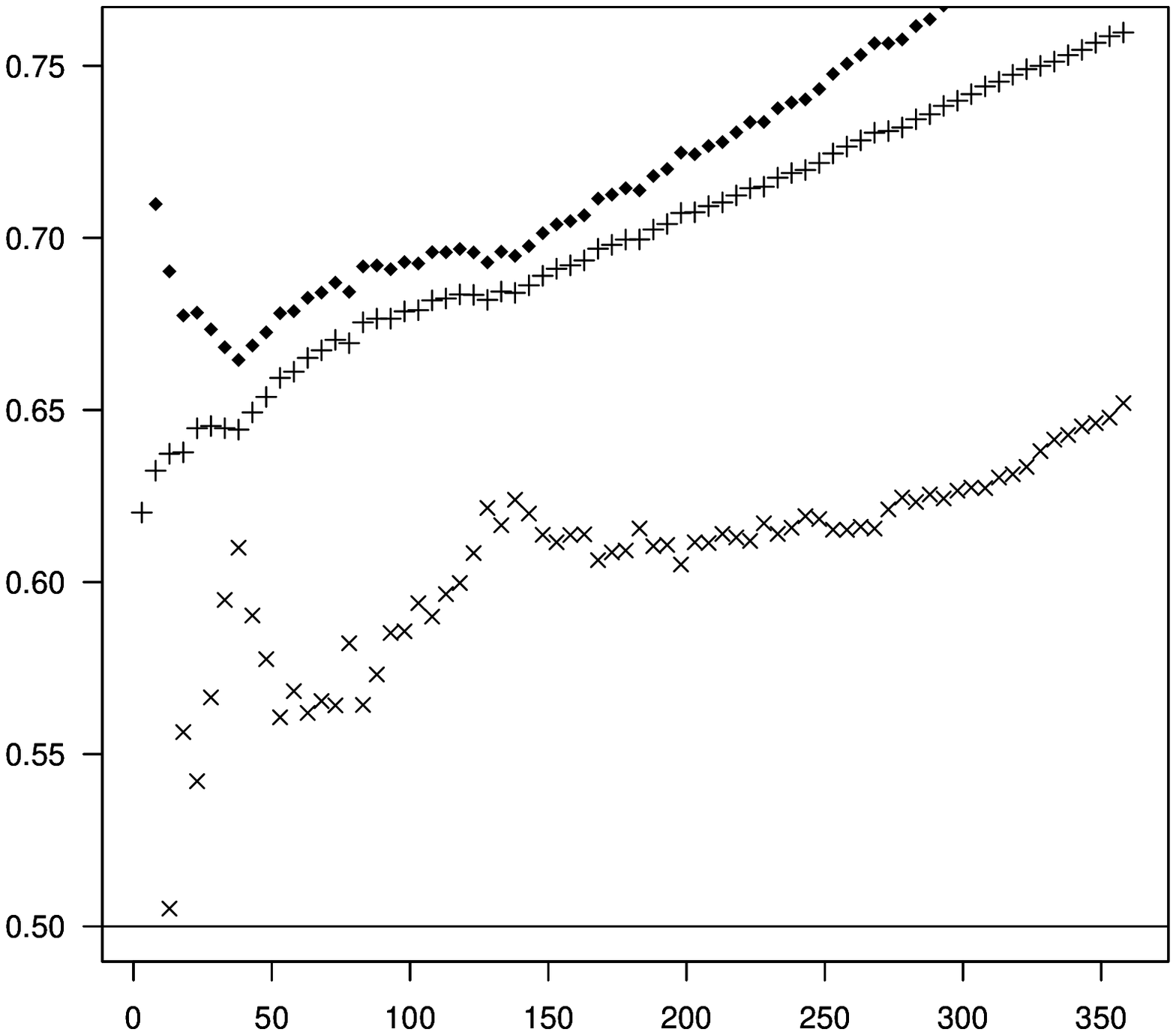}
\centerline{(a) Mean as a function of $k_n$}
\end{minipage}
\begin{minipage}{0.85\textwidth}
\includegraphics*[scale=0.5]{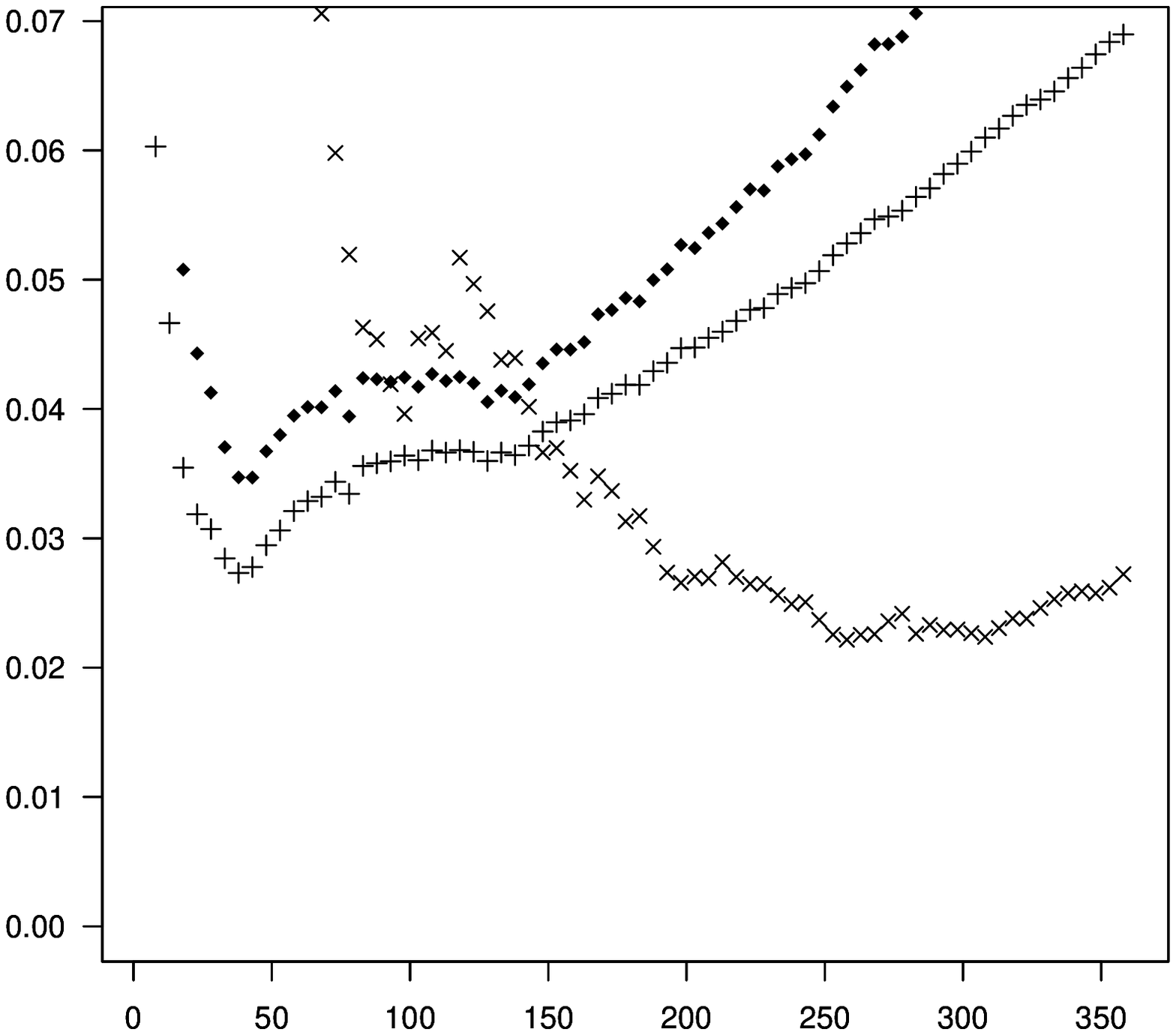}
\centerline{(b) Mean squared error as a function of $k_n$}
\end{minipage}
\end{center}
\caption{Comparison of estimates $\widehat{\theta}_n$ ($\times\times\times$), 
$\widetilde{\theta}_n$ ($\Diamondblack\Diamondblack\Diamondblack$) and $\check{\theta}_n$ ($+++$)  for the ${\mathcal{|N|}}(0,1)$ distribution. In (a), the straight line is the true value of $\theta$.}
\label{fignor}
\end{figure}

\begin{figure}
\begin{center}
\begin{minipage}{0.85\textwidth}
\includegraphics*[scale=0.5]{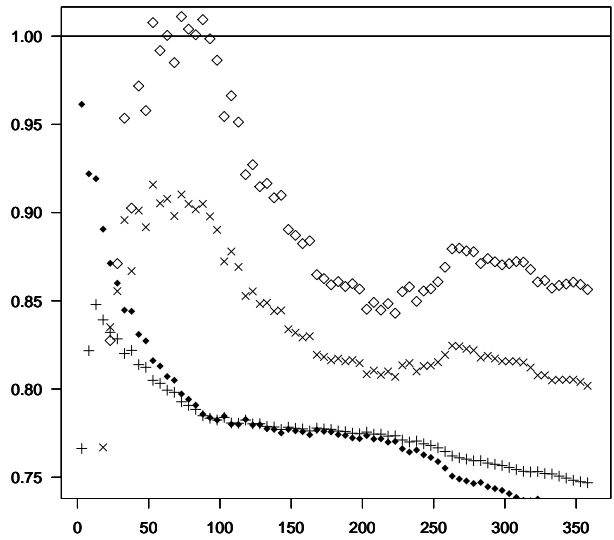}
\centerline{(a) Mean as a function of $k_n$}
\end{minipage}
\begin{minipage}{0.85\textwidth}
\includegraphics*[scale=0.5]{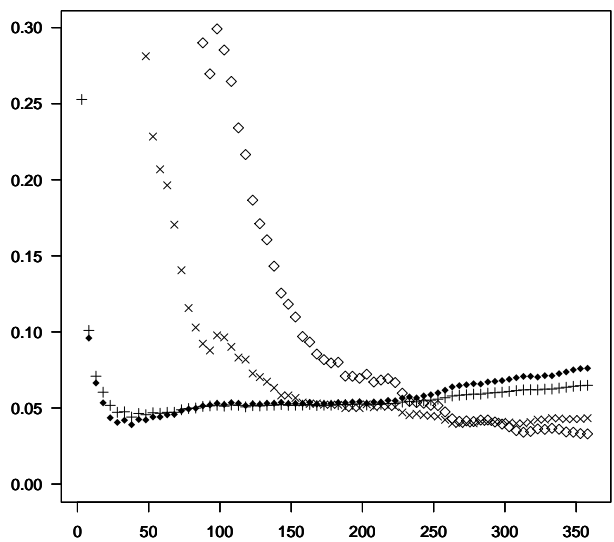}
\centerline{(b) Mean squared error as a function of $k_n$}
\end{minipage}
\end{center}
\caption{Comparison of estimates $\widehat{\theta}_n$
(with the canonical choice $\rho=-1$: $\times\times\times$), 
$\widehat{\theta}_n$ (with the true $\rho=-1/2$: 
$\Diamond\Diamond\Diamond$), 
$\widetilde{\theta}_n$ ($\Diamondblack\Diamondblack\Diamondblack$)
and $\check{\theta}_n$ ($+++$)  
for the ${\mathcal{D}}(1,0.5)$ distribution. In (a), the straight line is the true 
value of $\theta$.}
\label{fignew}
\end{figure}

\begin{figure}
\begin{center}
\begin{minipage}{0.85\textwidth}
\includegraphics*[scale=0.5]{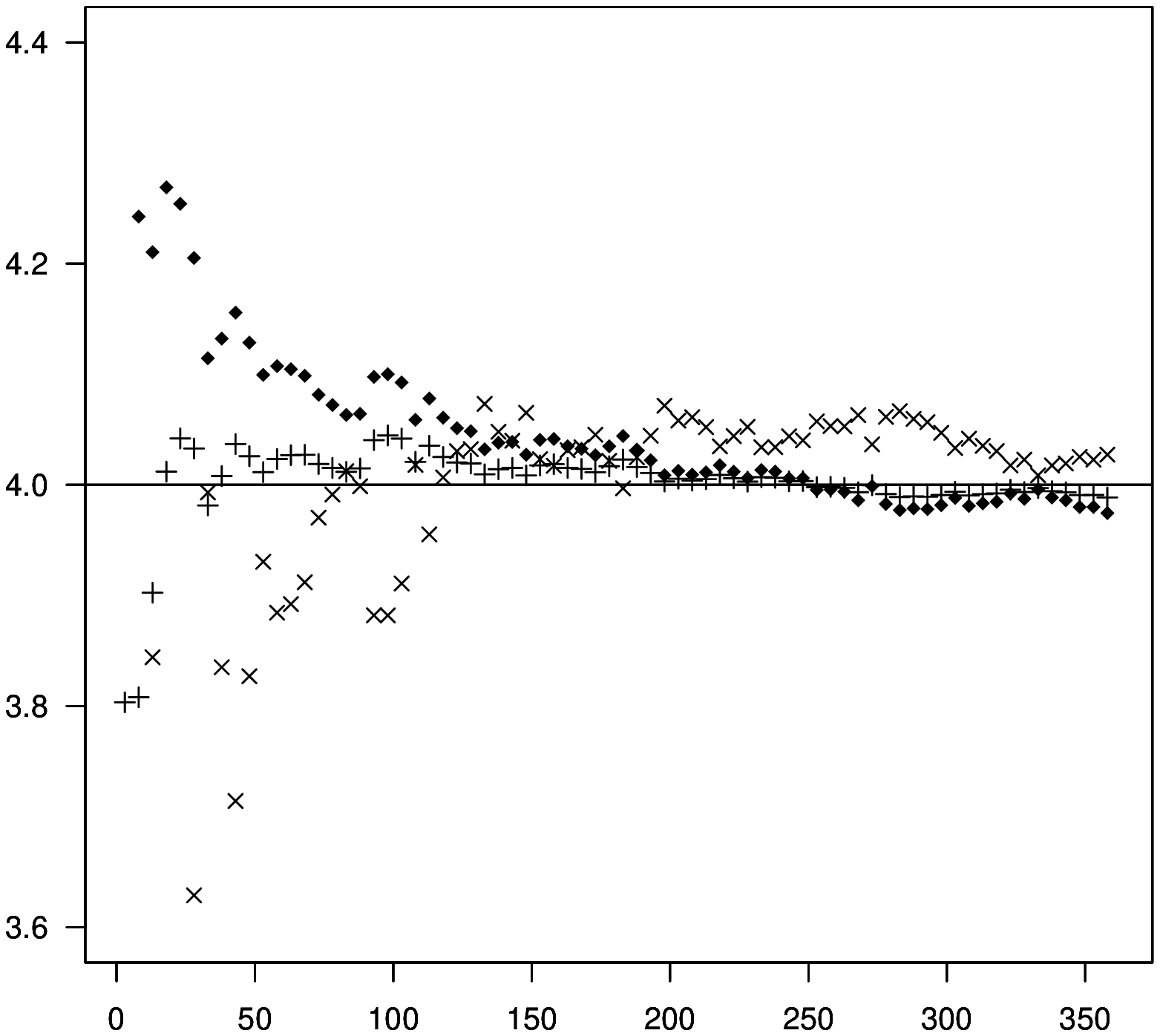}
\centerline{(a) Mean as a function of $k_n$}
\end{minipage}
\begin{minipage}{0.85\textwidth}
\includegraphics*[scale=0.5]{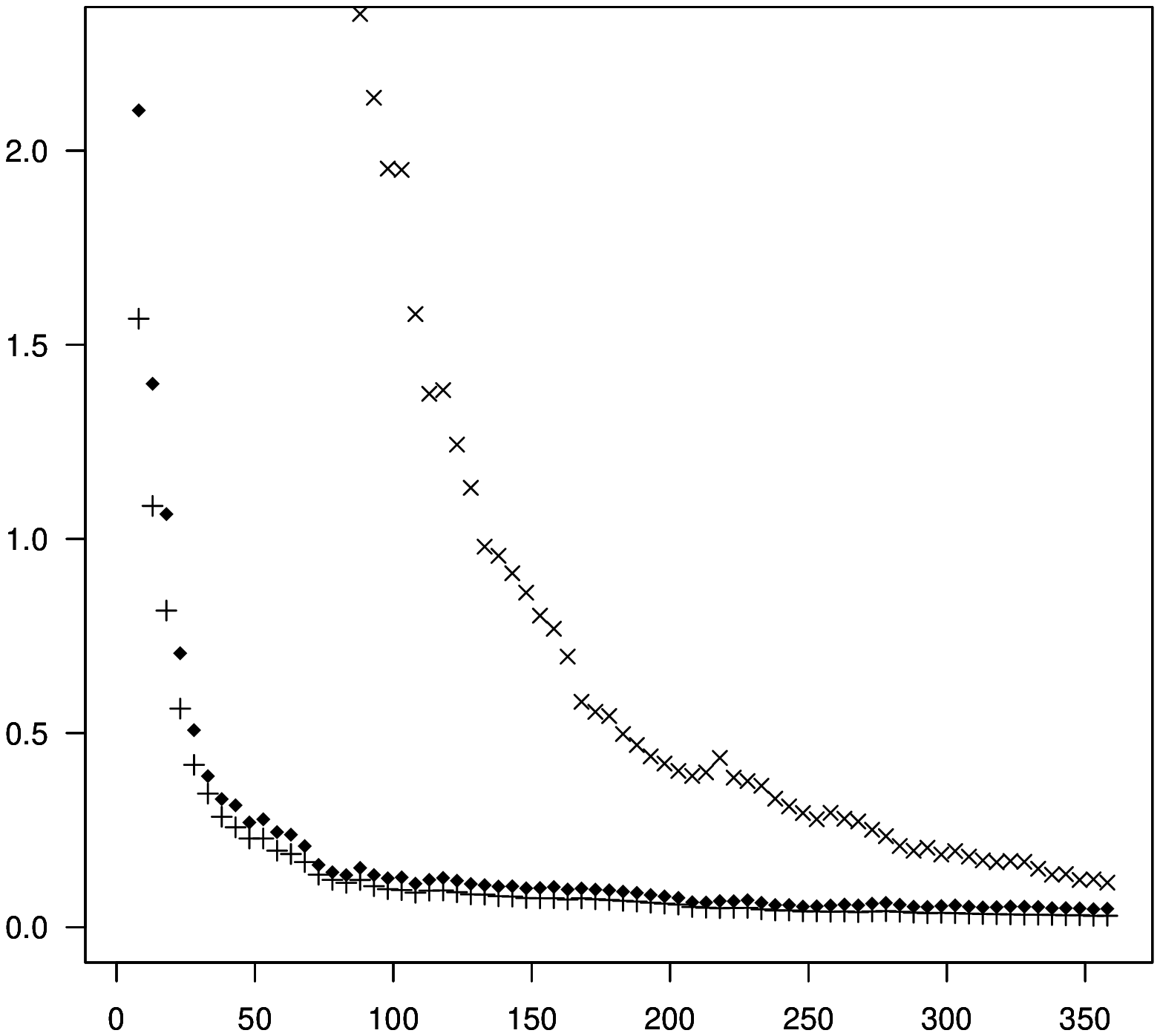}
\centerline{(b) Mean squared error as a function of $k_n$}
\end{minipage}
\end{center}
\caption{Comparison of estimates $\widehat{\theta}_n$ ($\times\times\times$), 
$\widetilde{\theta}_n$ ($\Diamondblack\Diamondblack\Diamondblack$) and $\check{\theta}_n$ ($+++$)  for the ${\mathcal{W}}(0.25,0.25)$ distribution. In (a), the straight line is the true value of $\theta$.}
\label{figweib025}
\end{figure}

\begin{figure}
\begin{center}
\begin{minipage}{0.85\textwidth}
\includegraphics*[scale=0.5]{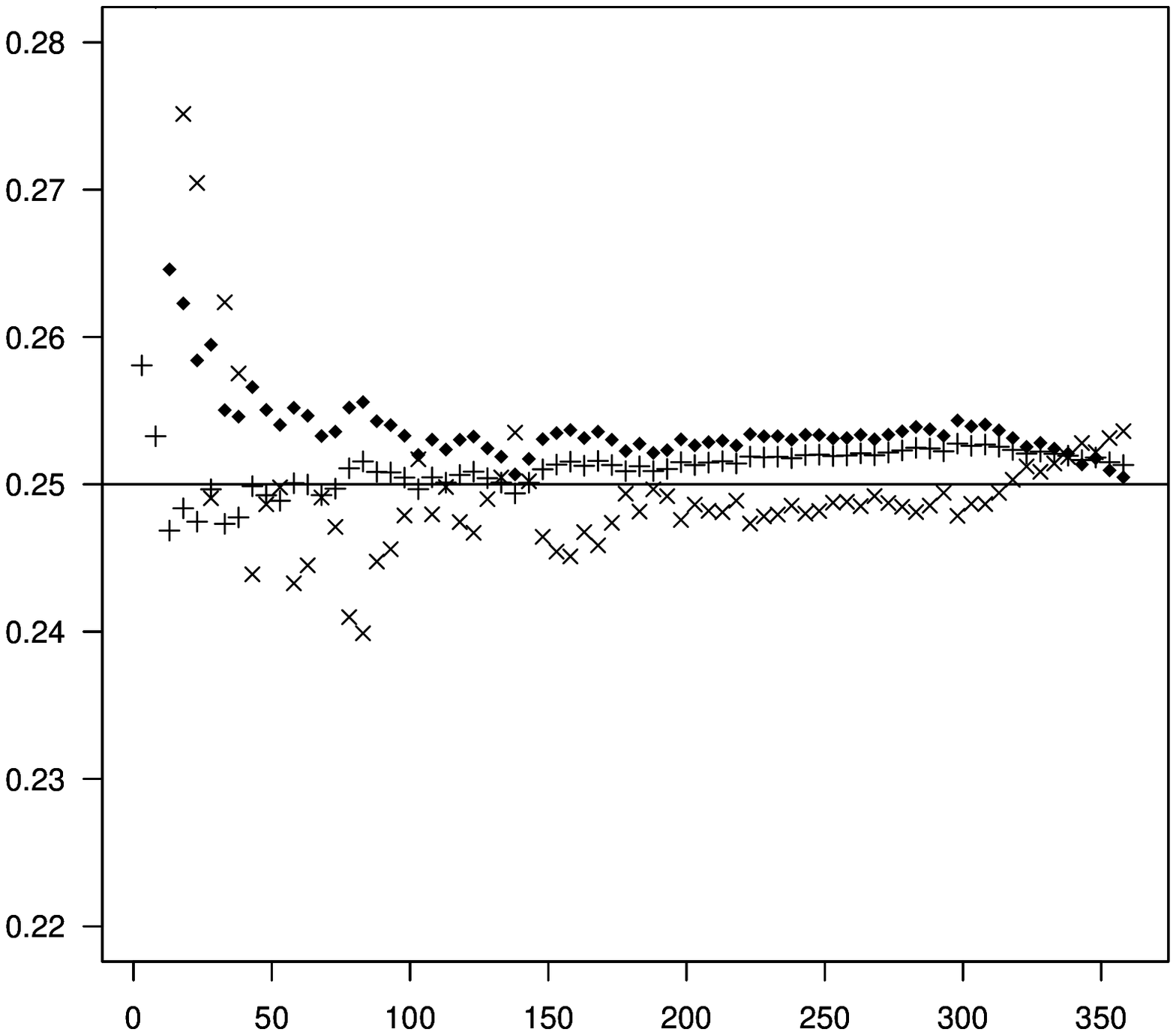}
\centerline{(a) Mean as a function of $k_n$}
\end{minipage}
\begin{minipage}{0.85\textwidth}
\includegraphics*[scale=0.5]{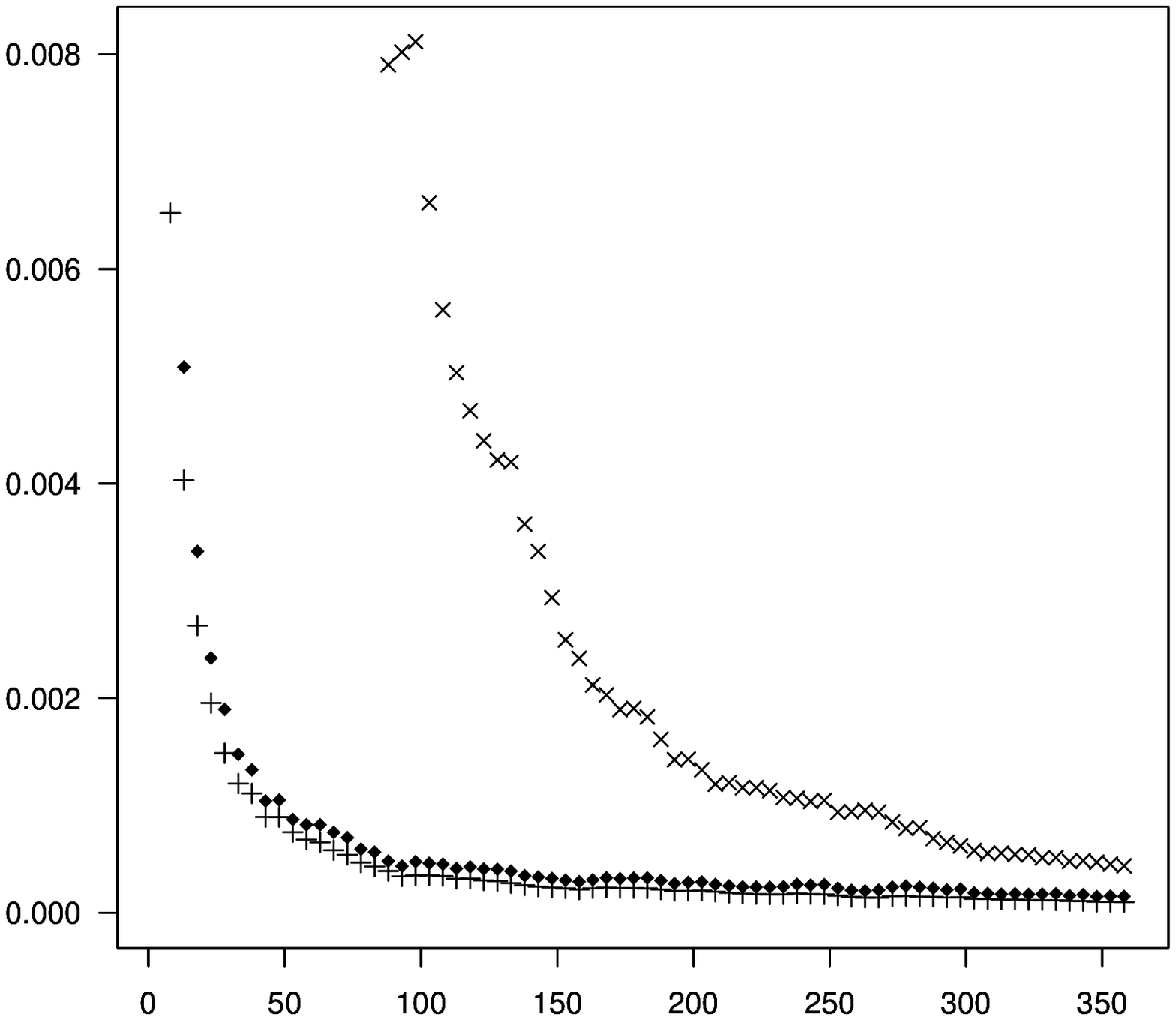}
\centerline{(b) Mean squared error as a function of $k_n$}
\end{minipage}
\end{center}
\caption{Comparison of estimates $\widehat{\theta}_n$ ($\times\times\times$), 
$\widetilde{\theta}_n$ ($\Diamondblack\Diamondblack\Diamondblack$) and $\check{\theta}_n$ ($+++$)  for the ${\mathcal{W}}(4,4)$ distribution. In (a), the straight line is the true value of $\theta$.}
\label{figweib4}
\end{figure}

\noindent The results are presented on figures~\ref{figgam025}--\ref{figweib4}. 
In all the plots, the graphs associated to ${\widetilde{\theta}}_n(k_n)$ and
${\check{\theta}}_n(k_n)$ are similar, with a slightly better behavior of ${\check{\theta}}_n(k_n)$.
The bias corrected estimator $\widehat{\theta}_n(k_n)$ always yields a smaller
bias than the two previous ones leading to better results for Gamma, Gaussian
and ${\cal D}$ distributions (figures~\ref{figgam025}--\ref{fignew}), 
even though a wrong value of $\rho$ is used (figure~\ref{fignew}).
On Weibull distributions, where the bias function is zero,
(figures~\ref{figweib025}--\ref{figweib4}), it presents a larger variance.  

\noindent Second, we investigate the behavior of the adaptive procedure
for selecting the number of upper order statistics in $\check \theta_n(k_n)$.
For $i=1,\dots,N$, we denote by
$$
 {\widehat{k_{n,i}}} = {\rm{arg}}\min_{k_n \in [1,350]} \widehat{AMSE}(\check 
 \theta_{n,i}(k_n))
$$
the value selected on the sample $({\mathcal X}_{n,i})$.
Note that, as in \cite{MB}, in our simulations,
we limited the range from which $k_n$ is selected to 
$\{1,\dots,350\}$. 
The mean and the standard deviation of this estimation
on the $N$ samples are given by
$$
\mu({{\widehat{k_n}}}) = \frac{1}{N} \sum_{i=1}^N {\widehat{k_{n,i}}}
\;\;\mbox{ and }\;\;\sigma({{\widehat{k_n}}}) = \sqrt{\frac{1}{N} \sum_{i=1}^N \left({\widehat{k_{n,i}}}-{\mu({{\widehat{k_n}}})}\right)^2}.
$$
As a comparison,
we introduce the value that would be obtained by
minimizing the true AMSE:
$$
 k_n^{opt}=  {\rm{arg}}\min_{k_n \in [1,350]} AMSE(\check \theta_{n}(k_n)).
$$
On each sample $({\mathcal X}_{n,i})$, the estimation of $\theta$
obtained with the selected parameter ${\widehat{k_{n,i}}}$ is
given by ${\check{\theta}}_{n,i}(\widehat{k_{n,i}})$.   The associated
empirical mean and standard deviation are:
\[ \mu({{{\check{\theta}}}_{n}}) = \frac{1}{N} \sum_{i=1}^N 
{\check{\theta}}_{n,i}(\widehat{k_{n,i}})\;\; \mbox{and} \;\;
\sigma({{\check{\theta_n}}})=\sqrt{\frac{1}{N} \sum_{i=1}^N 
\left({\check{\theta}}_{n,i}(\widehat{k_{n,i}})-\mu({{\check{\theta}}}_{n})\right)^2}. \]
Finally, to assess the quality of the selection procedure,
we compute the ratio $R_n$ of the empirical root mean squared error 
of ${\check{\theta}}_{n}(\widehat{k_{n}})$ and the minimal
empirical root mean squared error of ${\check{\theta}}_{n}({k_{n}})$:
\[ R_n^2 = \left.  {\displaystyle\sum_{i=1}^N 
({\check{\theta}}_{n,i}(\widehat{k_{n,i}})-\theta)^2} \right/{\displaystyle 
\min_{k_n \in [1,350]} \sum_{i=1}^N ({\check{\theta}}_{n,i}(k_n)-\theta)^2}. \]
Results are presented in Table~\ref{tabres}.
It appears that $R_n$ is usually ``close'' to 1, except for Weibull
distributions.  In this case, large values of $R_n$ together with
large values of $\mu(\hat k_n )$ indicate that the optimal
$k_n$ is larger than 350 observations.


\section{Real data}
\label{sec.cinq}

\begin{figure}
\begin{center}
\begin{minipage}{0.85\textwidth}
\includegraphics*[scale=0.5]{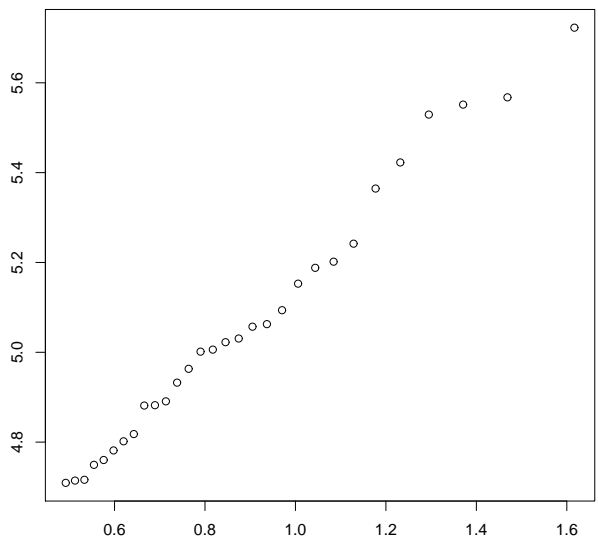}
\end{minipage}
\end{center}

\caption{Quantile-quantile plot obtained with $\widehat k_n=29$ on the Nidd river data.  } 
\label{nidd}
\end{figure}

Here, the good performance of the adaptive selection procedure 
is illustrated through the analysis of extreme events
on a benchmark real data set.
Nidd river data are
widely used in extreme value studies \cite{Hosking87,Dav}.
The raw data consist in $154$ exceedances of the level~$65$
m$^3$s$^{-1}$ by the river Nidd (Yorkshire, England)
during the period 1934-1969 (35 years).
The $N$-year return level is the water level which is exceeded
on average once in~$N$ years.
Hydrologists need to estimate extreme quantiles
in order to predict return levels over long periods.  
According to \cite{Hosking87}, the Nidd data may
reasonably be assumed to come from a distribution in the Gumbel
maximum domain of attraction.  This suggests to consider
Weibull tail-distributions as a possible model for such data.  
The adaptive selection procedure yields $\widehat k_n=29$.  
The resulting quantile-quantile plot (obtained by plotting the points
$(\log\log(n/i),\log{(X_{n-i+1,n})})$ for $i=1,\dots,\widehat k_n-1$)
is approximatively linear (see Figure~\ref{nidd}),
indicating a good fit of the Weibull tail-distribution 
for $x\geq X_{n-\hat k_n+1,n}$. 
We obtained
$\check \theta_n(\hat k_n) \simeq \widetilde \theta_n(\hat k_n) \simeq 0.89$.  
One can plug this result in the Weissman-type extreme quantile estimator 
proposed in \cite{Gar} to obtain $366m^3s^{-1}$
as an estimation of the 100-year return level.
Note that this result is in accordance with the results
obtained by profile-likelihood or Bayesian methods, see \cite{Dieb}
or \cite{Dav}.


\section{Concluding remarks}
\label{sec.six}

In this paper, we introduce a regression model, from which we derive a bias-reduced estimator for the Weibull tail-coefficient $\theta$.
Its asymptotic normality is established and an adaptive selection procedure
for $k_n$ is proposed.  
The efficiency of our approach is illustrated in a simulation study and
on a real data set. However, in many cases of practical interest, the problem of estimating a quantile $x_{p_n} = F^{-1}(1-p_n)$, with $p_n <1/n$, is much more important. Such a problem has already been studied in \cite{Gar} where the following Weissman-type estimator has been introduced
$$\widetilde x_{p_n}(k_n) = X_{n-k_n+1,n} \, \Biggr(\frac{\log (1 / p_n) }{ \log (n / k_n)}\Biggr)^{\widetilde \theta_n(k_n)}.$$ 
It is, however, desirable to refine $\widetilde x_{p_n}(k_n)$ with the additional information about the slowly varying function $\ell$ that is provided by the LS-estimates for $\theta$ and $b$. To this aim, condition $(R_{\ell}(b, \rho))$ is used to approximate the ratio $F^{-1}(1-p_n) / X_{n-k_n+1,n}$, noting that
$$X_{n-k_n+1,n} \stackrel{d}{=} F^{-1}(U_{n-k_n+1,n}),$$
with $U_{1,n} \leq ... \leq U_{n,n}$ the order statistics of a uniform $(0, 1)$ sample of size~$n$,
\begin{eqnarray*}
\frac{x_{p_n} }{ X_{n-k_n+1,n}} &\stackrel{d}{=} & \frac{F^{-1}(1-p_n) }{ F^{-1}(U_{n-k_n+1,n})} \\
&\stackrel{d}{=} & \frac{(-\log p_n)^{\theta} }{ (-\log(1-U_{n-k_n+1,n}))^{\theta}} \, \frac{\ell(-\log p_n) }{ \ell(-\log(1-U_{n-k_n+1,n}))}\\
&\stackrel{d}{\simeq} & \Biggr(\frac{\log (1 / p_n) }{ \log (n / k_n)}\Biggr)^{\theta} \, \exp\Biggr[b(\log (n / k_n)) \frac{\Bigr(\frac{\log (1 / p_n) }{ \log (n/k_n)}\Bigr)^{\rho} - 1 }{ \rho}\Biggr].
\end{eqnarray*}
The last step follows from replacing $U_{k_n+1,n}$ (resp. $E_{n-k_n+1,n}$) by $k_n/n$ (resp. $\log(n/k_n)$). Hence, we arrive at the following estimator for extreme quantiles
$$\widehat x_{p_n}(k_n) = X_{n-k_n+1,n} \, \Biggr(\frac{\log (1 / p_n) }{ \log (n/k_n)}\Biggr)^{\widehat \theta_n(k_n)} \, \exp\Biggr[\widehat b(\log (n / k_n)) \frac{\Bigr(\frac{\log (1 / p_n) }{ \log (n / k_n)}\Bigr)^{\widehat \rho} - 1 }{ \widehat \rho}\Biggr],$$
which is similar to the estimator proposed by \cite{MB2}
in the case of Pareto-type distributions.  
Here, the LS-estimators of $\theta$ and $b$ can be used after substituting 
$\rho$ by the canonical choice $-1$. The study of the asymptotic properties of 
such an estimator is the aim of \cite{Dieb2}.


\section{Proofs of our results}\label{sec.sept}

\subsection{Preliminary lemmas}

\begin{lemma}
\label{lemEnj}
For all $1\leq j \leq k_n$ such that $k_n \to \infty$ and ${k_n / n}\to 0$, we have
\begin{eqnarray*}
\frac{E_{n-j,n} }{ \log(n / j)} = 1 +  O_{\mathbb P}\Biggr(\frac{1}{ \log(n / k_n)}\Biggr)\quad \mbox{{\it uniformly in}}\,  j.
\end{eqnarray*}
\end{lemma}
\noindent
{\bf Proof of Lemma \ref{lemEnj}.} According to R\'enyi's representation, we have
$$E_{n-j,n} \stackrel{d}{=} \sum_{\ell=1}^{n-j+1} \frac{f_{n-\ell-j+1} }{ \ell+j-1}$$
where $f_j  \stackrel{i.i.d.}{\sim} \mbox{Exp}(1)$. Since 
$$\mbox{Var}\Biggr(\sum_{\ell=1}^{n-j+1} \frac{f_{n-\ell-j+1} }{ \ell+j-1}\Biggr) = O(1),$$
denoting 
$$T_{j,n} := \sum_{\ell=1}^{n-j+1} \Biggr[\frac{f_{n-\ell-j+1} }{ \ell+j-1} - \mathbb E\frac{f_{n-\ell-j+1} }{ \ell+j-1}\Biggr],$$
we have, using Kolmogorov's inequality \cite{Shor} (p.183), that
$$\mathbb P\left(\max_{1\leq j\leq k_n} \vert T_{j,n} \vert \geq \lambda\right) \leq \frac{\mbox{Var} (T_{1,n}) }{ \lambda^2}, \qquad \lambda>0.$$
This implies that $T_{j,n} = O_{\mathbb P}(1)$ uniformly in $j$. Taking into account the fact that
$$\Biggr\vert \sum_{\ell=j}^n \frac{1}{ \ell} - \log (n / j)\Biggr\vert = O(1) \qquad \mbox{uniformly in} \, j, 1\leq j \leq k_n,$$
it is easy to deduce Lemma \ref{lemEnj}. \hfill{\cqfd}

\noindent Let us introduce the $E_m-$function defined by the integral
$$E_m(x) := \int_1^{\infty} \frac{e^{-xt}}{ t^m} dt$$
for a positive integer $m$. The asymptotic expansion of this integral is given in the following lemma.
\begin{lemma}
\label{lemEn}
As $x\to \infty$, for any fixed positive integers $m$ and $p$, we have
$$E_m(x) = \frac{e^{-x} }{ x} \Biggr\{1-\frac{m }{ x} +  ... + (-1)^p \frac{m(m+1)...(m+p-1)}{ x^p} + O\Bigr(\frac{1 }{ x^{p+1}}\Bigr)\Biggr\}.$$
\end{lemma}
\noindent
The proof of this lemma is straightforward from \cite{Abra} p. 227-233 and the $O-$term can be obtained by a Taylor expansion with an integral remainder. 
Denote
$$\mu_p := \frac{1 }{ k_n} \sum_{j=1}^{k_n} \left(x_j-\overline x_{k_n}\right)^p, \, p \in {\mathbb N}^*.$$
The next lemma provides a first order expansion of this Riemman sum.
\begin{lemma}
\label{lemmup}
If $k_n \to \infty$, ${k_n / n}\to 0$, $\frac{k_n }{ \log(n / k_n)} \to \infty$ and $\frac{\log^2 k_n }{  \log(n / k_n)} \to 0,$ then
$$\mu_p \sim C_p ({\log(n / k_n}))^{-p} \mbox{ as $n\to \infty$, where } C_p = \int_0^1(\log x+1)^p dx < \infty.$$
\end{lemma}
\noindent
{\bf Proof of Lemma \ref{lemmup}.} Denote $\alpha_n = \frac{1 }{ \log(n/k_n)}$.
 Then $\overline x_{k_n}$ can be rewritten as 
$$\overline x_{k_n}= \frac{1}{k_n}+ \Biggr(\frac{1}{ k_n} \sum_{j=1}^{k_n-1} f_n(j/k_n)-\int_0^1 f_n(x)dx\Biggr)+\int_0^1 f_n(x)dx =: \frac{1}{k_n}+ T_1 + T_2,$$
where $f_n(x) = (1-\alpha_n \log x)^{-1}, x \in [0, 1].$
Denoting by $f_n^{(i)}$, $i\in\{1,2\},$ the $i$th derivative of $f_n$,
we infer that
\begin{eqnarray*}
T_1& =& \sum_{j=1}^{k_n-1} \int_{j/k_n}^{(j+1)/k_n} ({j / k_n}-t) f_n^{(1)}({j/ k_n})dt \\ &+&
\sum_{j=1}^{k_n-1} \int_{j/k_n}^{(j+1)/k_n} \int_{j/k_n}^t (x-t) f_n^{(2)}(x) dxdt+ \int_0^{1/k_n} f_n(x)dx\\
&=:& T_3+T_4+T_5.
\end{eqnarray*}
Remark that
\begin{eqnarray*}
T_3 & = & -\frac{1 }{ 2 k_n} \Biggr ( \frac{1 }{ k_n} \sum_{j=1}^{k_n-1} f_n^{(1)}
 ( {j / k_n} )- \int_{1/k_n}^1 f^{(1)}_n(t)dt \Biggr ) -\frac{1 }{ 2 k_n} \int_{1/k_n}^1 f_n^{(1)}(t) dt \\
 \ & =: & - \frac{1 }{ 2 k_n} T_6 + T_7.
\end{eqnarray*}
Since $f_n^{(1)}$ is positive and decreasing on $\Bigr[\frac{1}{ k_n}, 1\Bigr]$ for $n$ sufficiently large, we can prove that
\begin{eqnarray*}
\vert T_4\vert  &\leq&  \frac{1}{ 2 k^2_n} \Bigr\vert f_n^{(1)}({1/ k_n}) - f_n^{(1)}(1)\Bigr\vert = o({1/ k_n}), \\
T_5 &=& O({1/ k_n}),\\ 
\vert T_6 \vert &\leq& \frac{1 }{ k_n} \Bigr \vert f_n^{(1)}({1/ k_n}) - f_n^{(1)}(1)\Bigr\vert =o(1),\\ 
T_7 &=& -\frac{1}{ 2 k_n} \Bigr ( f_n(1)-f_n  ( {1 / k_n}  )\Bigr ) = o({1/ k_n}),
\end{eqnarray*}
and consequently $T_1 = O({1/ k_n}).$ 
Besides, a direct application of Lemma~\ref{lemEn} provides
$$T_2 = 1-\alpha_n+O(\alpha_n^2).$$
Therefore $\overline x_{k_n} = 1-\alpha_n + O({1 / k_n})+O(\alpha_n^2).$ Now, we can check that
$$\mu_p = \alpha_n^p \Biggr\{\frac{1}{ k_n} \sum_{j=1}^{k_n} (\log(j / k_n))  + 1)^p+R_n\Biggr\}$$
where 
$$R_n = \frac{1}{ k_n} \sum_{j=1}^{k_n-1} \Biggr\{(\log({j / k_n}) + 1 + \varepsilon_n)^p - \Bigr(\log({j / k_n})  + 1 \Bigr)^p\Biggr\}$$
with $\varepsilon_n = O\left(\alpha_n \log^2 k_n\right)+O\left(\frac{1}{ k\alpha_n}\right)$ which tends to 0 by assumption.
Since $\frac{1 }{ C_p} \frac{1}{ k_n} \sum_{j=1}^{k_n} (\log({j / k_n}) + 1 )^p \to 1$, in order to conclude  the proof of Lemma~\ref{lemmup}, we only have to remark that $R_n\to 0$. \hfill{\cqfd}

\subsection{ Proof of Lemma \ref{lemAn}}

\noindent Remark that
\begin{eqnarray*}
\alpha_{j,n}&:=& j \log(n/ j) \Bigr[A_{n}(j)-A_{n}(j+1)\Bigr] \\
&=& \theta j \log(n/j) \log (E_{n-j+1,n}/ E_{n-j,n})\\
&=& \theta \log(n/ j) \, {j(E_{n-j+1,n}- E_{n-j,n}) / E^*_{n-j,n}}\\
&\stackrel{d}{=}& \theta f_j {\log(n/ j) / E^*_{n-j,n}}
\end{eqnarray*}
\noindent
where $E^*_{n-j,n} \in [E_{n-j,n}; E_{n-j+1,n}].$ 
Consequently, from Lemma~\ref{lemEnj},
\begin{eqnarray}
\nonumber
\alpha_{j,n}&=& \theta f_j + O_{\mathbb P}\Biggr(\frac{1}{ \log(n / k_n)}\Biggr)\\
&=& \theta f_j + o_{\mathbb P}\Bigr(b(\log(n /k_n))\Bigr),
\label{Q1}
\end{eqnarray}
\noindent
by the assumption $x|b(x)| \to \infty$ as $x\to \infty$ with a $o_{\mathbb P}-$term which is uniform in $j$. Lemma~\ref{lemAn} is therefore proved. \hfill{\cqfd}

\subsection{Proof of Lemma \ref{lemBn}}

\noindent We consider
$$\beta_{j,n} :=  j \log(n/ j) \Bigr[B_{n}(j)-B_{n}(j+1)\Bigr].$$
\noindent
In order to study this term, we will use the notations 
$\lambda_{1j} = {E_{n-j+1, n} / E_{n-k_n+1,n}}$, 
$\lambda_{2j} = {E_{n-j, n} / E_{n-k_n+1,n}}$ and $y_{k_n} = E_{n-k_n+1,n}$, and we rewrite $\beta_{j,n}$ as
$$
\beta_{j,n} = j \, \log(n / j) \Biggr\{\log \ell\Bigr(\lambda_{2j} \, \frac{\lambda_{1j} }{ \lambda_{2j}} \, y_{k_n}\Bigr) - \log \ell\Bigr(\lambda_{2j} \, y_{k_n}\Bigr)\Biggr\}.$$
\noindent
It is clear that $1\leq {\lambda_{1j} / \lambda_{2j}}  \stackrel{\mathbb P}{\longrightarrow} 1$ uniformly in $j$ by Lemma~\ref{lemEnj} and therefore for $n \geq N_0$, ${\lambda_{1j} / \lambda_{2j}} \in [1, 2]$ in probability. Under our assumption $(R_{\ell}(b, \rho))$ on the slowly varying function, we deduce that
$$
\beta_{j,n} = j \, \log(n / j) \Bigr\{b(\lambda_{2j}y_{k_n}) K_{\rho}({\lambda_{1j} / \lambda_{2j}}) (1+o_{\mathbb P}(1))\Bigr\}.
$$
Now, since $\lambda_{2j}\stackrel{\mathbb P}{\longrightarrow} 1$ uniformly in $j$ and $b(.)$ is regularly varying with index $\rho$,  $b(\lambda_{2j} y_{k_n}) = \lambda_{2j}^{\rho} b(y_{k_n}) (1+o_{\mathbb P}(1))$ with a $o_{\mathbb P}(1)$-term uniform in $j$.
\noindent
Therefore
$$\beta_{j,n} = j \, \log(n / j) \,  b(y_{k_n}) \Bigr\{ \lambda_{2j}^{\rho}\, K_{\rho}(  {\lambda_{1j} / \lambda_{2j}}) (1+o_{\mathbb P}(1)) \Bigr\}.$$
\noindent
Again, uniformly in $j$,
$$K_{\rho}(  {\lambda_{1j} / \lambda_{2j}}) = \left({\lambda_{1j}/ \lambda_{2j}}-1\right)(1+o_{\mathbb P}(1)),$$
\noindent
which implies that $\beta_{j,n}$ can be rewritten as follows:
\begin{eqnarray*}
\beta_{j,n} =  -\, j \, \log(n / j) \,  b(y_{k_n}) (\lambda_{2j}- \lambda_{1j}) \lambda_{2j}^{\rho-1}(1+o_{\mathbb P}(1)).
\end{eqnarray*}
\noindent
Therefore, we have
$$\beta_{j,n} = f_j \Biggr(\frac{\log (n / j)}{ \log(n /k_n)}\Biggr)^{\rho} b(y_{k_n}) (1+o_{\mathbb P}(1)),$$
with a $o_{\mathbb P}(1)$-term which is uniform in $j$.
This achieves the proof of Lemma~\ref{lemBn}. \hfill{\cqfd}

\noindent
Remark that, since $\frac{\log(n/j)}{\log(n/k_n)}\to 1$ uniformly in $j$,
one also has 
$$\beta_{j,n} = f_j \Biggr(\frac{\log (n / j)}{ \log(n / k_n)}\Biggr)^{-1} b(y_{k_n}) (1+o_{\mathbb P}(1)),$$
with a $o_{\mathbb P}(1)$-term which is uniform in $j$,
and this proves Corollary~\ref{Coro}.

\subsection{Proof of Theorem~\ref{thdeux}}

\noindent
From model (\ref{regression}), we infer that
\begin{eqnarray*}
&&\sqrt{k_n}\, \Biggr(\check{\theta}_n(k_n) - \theta-b(\log (n /k_n)) \frac{1 }{ k_n} \sum_{j=1}^{k_n} \left(\frac{\log (n / j) }{ \log (n /k_n)}\right)^{\rho}\Biggr)\\ 
&&= \sqrt{k_n} \, \theta \, \frac{1 }{ k_n} \sum_{j=1}^{k_n} (f_j-1) + \sqrt{k_n} b(\log (n / k_n)) \frac{1 }{ k_n} \sum_{j=1}^{k_n} \left(\frac{\log (n / j) }{ \log (n / k_n)}\right)^{\rho} (f_j-1)\\ 
&&+ o_{\mathbb P}\left(\sqrt{k_n}\, b(\log (n / k_n))\right).
\end{eqnarray*}
Now, an application of Tchebychev's inequality gives that 
$$\frac{1 }{ k_n} \sum_{j=1}^{k_n} \left(\frac{\log (n / j) }{ \log (n / k_n)}\right)^{\rho} (f_j-1)=o_{\mathbb P}(1).$$
Then, under our assumptions, Theorem~\ref{thdeux} follows by an application of the Central Limit Theorem.\hfill{\cqfd}

\subsection{Proof of Theorem~\ref{thtrois}}

From Corollary~\ref{Coro}, we have
\begin{eqnarray*}
&&\frac{\sqrt{k_n} }{ \log(n / k_n)} \Bigr(\widehat \theta_n(k_n)-\theta\Bigr)\\
&&= \frac{\sqrt{k_n} }{ \log(n /k_n)} \frac{1 }{ k_n} \sum_{j=1}^{k_n} \left(\theta+b(\log(n / k_n))x_j\right)\left(1-\frac{x_j-\overline x_{k_n} }{ \mu_2}\overline x_{k_n}\right)(f_j-1)\\
& &+ o_{\mathbb P}\left(\frac{\sqrt{k_n} }{ \log(n / k_n)} b(\log(n /k_n)\Bigr)\right).
\end{eqnarray*}
\noindent
Since we have (\ref{cond2}) and (\ref{cond2bis}), the $o_{\mathbb P}$-term is negligible. The first term can be viewed as a sum of a weighted mean of independent and identically distributed variables. 
Now, using Lyapounov's theorem, we only have to show that
$$\lim_{k_n \to \infty} \frac{1 }{ s_{k_n}^4} \sum_{j=1}^{k_n} \mathbb E X_j^4 = 0,$$
\noindent
where $X_j = \Bigr(\theta+b(\log(n / k_n)) x_j\Bigr) \Bigr(1-\frac{x_j-\overline x_{k_n} }{ \mu_2} \overline x_{k_n}\Bigr)(f_j-1)$, $j=1, ..., k_n$ and $s_{k_n}^2 = \sum_{j=1}^{k_n} \mbox{Var} X_j.$
\noindent
We remark that 
$$s_{k_n}^2 \sim \theta^2 \sum_{j=1}^{k_n} \left(1-\frac{x_j-\overline x_{k_n} }{ \mu_2} \overline x_{k_n}\right)^2\quad \mbox{as $n\to\infty$}$$
\noindent
and
$$\sum_{j=1}^{k_n} \mathbb E X_j^4 \sim 9 \theta^4 \sum_{j=1}^{k_n} \left(1-\frac{x_j-\overline x_{k_n} }{ \mu_2} \overline x_{k_n}\right)^4\quad \mbox{as $n\to\infty$}$$
from which we deduce by direct computations that
\begin{eqnarray*}
\frac{1 }{ s_{k_n}^4} \sum_{j=1}^{k_n} \mathbb E X_j^4 &\sim& \frac{9 }{ k_n} \, \frac{\mu_2^4+6(\overline x_{k_n})^2\mu_2^3-4(\overline x_{k_n})^3\mu_2\mu_3 + (\overline x_{k_n})^4\mu_4 }{ [\mu_2^2+(\overline x_{k_n})^2\mu_2]^2}\\
&\sim& \frac{9 C_4 }{ k_n}
\end{eqnarray*}
by Lemma~\ref{lemmup}.
Our Theorem~\ref{thtrois} now follows from the fact that
$$s_{k_n}^2 \sim \theta^2 k_n \log^2(n/k_n).$$
\hfill{\cqfd}


\section*{Appendix}

In this appendix, we briefly show how to adapt
Hall's class of distribution function \cite{Hall} to the
framework of Weibull tail-distributions.  
We introduce the class of distributions ${\cal D}(\alpha,\beta)$
with distribution function given by
$$
1-F(x) = \exp(-H(x))  \mbox{ where }  H^{-1}(x):= 
x^{1/\alpha} (1+x^{-\beta}),
$$
$\alpha$ and $\beta$ being two parameters such that
\begin{equation}
\label{condD}
0<\alpha,\; 0< \beta< 1\mbox{ and } \alpha\beta\leq 1.
\end{equation}
It is easily seen that under (\ref{condD}), the above
class of distributions fulfill assumptions
(\ref{model}) with $(R_{\ell}(b, \rho))$ and (\ref{hypob})
where $\theta=1/\alpha$, $\rho=-\beta$,
$\ell(x)=1+x^{-\beta}$ and
$b(x)=-\beta x^{-\beta}$.  It is thus possible
to obtain distributions with arbitrary $\theta>0$
and $-1<\rho<0$.   These results are
summarized in Table~\ref{tabex}.

\section*{Acknowledgement}

The authors are very grateful to the referees for a careful reading
of the paper that led so significant improvements of the earlier
draft.  


\newpage


\begin{table}
\begin{center}
$
\begin{array}{|c|c|c|c|}
\hline
&&&\\
\mbox{Distribution}                       & \theta & b(x) & \rho \\
&&&\\
\hline
&&&\\
\mbox{Absolute Gaussian }|{\mathcal N}|(\mu,\sigma^2)      & 1/2   & \displaystyle\frac{1}{4} \frac{\log x}{x} & -1 \\
&&&\\
\mbox{Gamma }\Gamma(\alpha\neq 1,\beta)      & 1     & (1-\alpha) \displaystyle \frac{\log x}{x} & -1 \\
&&&\\
\mbox{Weibull }{\mathcal W}(\alpha,\lambda)    & 1/\alpha & 0 & -\infty \\
&&&\\
{\mathcal D}(\alpha,\beta)    & 1/\alpha & -\beta x^{-\beta} & -\beta \\
&&&\\
\hline
\end{array}
$
\end{center}
\caption{Parameters $\theta$, $\rho$ and the function $b(x)$ associated
to some  distributions}
\label{tabex}
\end{table}

\begin{table}

\begin{center}
\begin{tabular}{|c | c | c | c | c | c | c | c | c |}
\hline 
 \ & \ & \ & \ & \ & \ & \ & \ & \ \\
Distribution & $\theta$ & $\rho$ & ${\mu({{\widehat{k_n}}})}$ & $\sigma({{\widehat{k_n}}})$ & $\mu({{\check{\theta}}}_{n})$ & $\sigma({\check{\theta_n}})$ & $R_n$ & $k_n^{opt}$ \\
 \ & \ & \ & \ & \ & \ & \ & \ & \ \\
\hline 
 \ & \ & \ & \ & \ & \ & \ & \ & \ \\
$\Gamma (0.25,1)$ & 1 & -1 & 105.5 & 62.2 & 1.667 & 0.294 & 1.26 & 186 \\
$\Gamma (4,1)$ & 1 & -1 & 222.7 & 82.1 & 0.548 & 0.051 & 1.13 & 184 \\
${|{\cal N}|}(0,1)$ & 0.5 & -1 & 246.6 & 81.1 & 0.679 & 0.109 & 1.21 & 189 \\
${{\cal W}}(0.25,0.25)$ & 4 & -$\infty$ & 305.8 & 59.0 & 4.016 & 0.265 & 1.62 & 350 \\
${{\cal W}}(4,4)$ & 0.25 & -$\infty$ & 310.4 & 50.9 & 0.249 & 0.013 & 1.43 & 350 \\
${{\cal D}(1,0.5)}$ & 1 & -0.5 & 281.5 & 71.1 & 0.789 & 0.053 & 1.14 & 43 \\
 \ & \ & \ & \ & \ & \ & \ & \ & \ \\
\hline 
\end{tabular}
\end{center}
\caption{Simulation results of the adaptive selection procedure}
\label{tabres}
\end{table}

\end{document}